\documentclass[12pt]{article}

\makeatletter
 
 \@addtoreset{equation}{section}
\makeatother

\newcommand{\ul}{\underline}

\newcommand{\prf}{\noindent{\bf Proof. }}

\newcommand{\rem}{\noindent{\bf Remark. }}

\newcommand{\ff}[1]{%
    {\bf F}_{#1}}

\newcommand{\qed}{\hbox{\rule[-2pt]{5pt}{11pt}}}

\newtheorem{dfn}{Definition}[section]
\newtheorem{thm}[dfn]{Theorem}
\newtheorem{prop}[dfn]{Proposition}

\newtheorem{lem}[dfn]{Lemma}
\newtheorem{exam}[dfn]{Example}

\newtheorem{prob}[dfn]{Problem}

\begin{document}
\title{Zeta functions for formal weight enumerators and an analogue of the Mallows-Sloane bound}
\author{Koji Chinen\footnotemark[1]}
\date{}
\maketitle
\begin{abstract}
In 1999, Iwan Duursma defined the zeta function for a linear code as a generating function of its Hamming weight enumerator. It has various properties similar to those of the zeta function of an algebraic curve. 

This article extends Duursma's theory to the case of the formal weight enumerators and shows that there exists a similar structure to that of the weight enumerators of type II codes. The extremal property of the formal weight enumerators is also considered and an analogue of the Mallows-Sloane bound is deduced. 
\end{abstract}
\footnotetext[1]{Department of Mathematics, Faculty of Engineering, Osaka Institute of Technology. 5-16-1 Omiya, Asahi-ku, Osaka 535-8585, Japan. E-mail: YHK03302@nifty.ne.jp}



\section{Introduction}
Let $p$ be a prime, $q=p^r$ for some positive integer $r$ and we denote the finite field with $q$ elements by $\ff q$. Let $C$ be an $[n,k,d]$-code over $\ff q$ with the Hamming weight enumerator $W_C(x,y)$. In 1999, Duursma \cite{Du1} defined the zeta function for a linear code as follows:
\begin{dfn}[Duursma]\label{dfn:zeta}
For any linear code $C$, there exists a unique polynomial $P(T)\in{\bf Q}[T]$ of degree at most $n-d$ such that
$$\frac{P(T)}{(1-T)(1-qT)}(y(1-T)+xT)^n=\cdots +\frac{W_C(x,y)-x^n}{q-1}T^{n-d}+ \cdots.$$
We call $P(T)$ the zeta polynomial of the code $C$, and $Z(T)=P(T)/((1-T)(1-qT))$ the zeta function of $C$. 
\end{dfn}
In his subsequent papers \cite{Du2, Du3, Du4}, Duursma deduces various interesting properties of $P(T)$ and discusses their possible applications to the coding theory. Among them, the functional equation and a Riemann hypothesis analogue for self-dual codes attract interests of many mathematicians, both in coding theory and number theory. When $C$ is self-dual, the MacWilliams identity leads the functional equation of $P(T)$ of the form
\begin{equation}\label{eq:func_eq}
P(T)=P\Bigl(\frac{1}{qT}\Bigr)q^g T^{2g}
\end{equation}
($g=n+1-k-d$, see \cite[p.59]{Du2}) and we can formulate an analogue of the  Riemann hypothesis as follows (see Duursma \cite[Definition 4.1]{Du3}): 
\begin{dfn}\label{dfn:RH}
The code $C$ satisfies the Riemann hypothesis if all the zeros of $P(T)$ have the same absolute value $1/\sqrt q$. 
\end{dfn}
One of the striking differences between the zeta functions of self-dual codes and those of algebraic curves is that the Riemann hypothesis for self-dual codes often fails to hold (in the case of the algebraic curves, the Riemann hypothesis is always true, as was proved by A. Weil). Finding the equivalent condition for a self-dual code to satisfy the Riemann hypothesis is still an open problem, but Duursma proposes a certain sufficient condition (see \cite[Open Problem  4.2]{Du3}): 
\begin{prob}\label{prob:RH}
Prove or disprove that all extremal weight enumerators satisfy the Riemann hypothesis. 
\end{prob}
A self-dual code $C$ is called extremal if it has the largest possible minimum distance (see Pless \cite[p.139]{Pl}). There are 4 sequences of extremal self-dual codes (Types I, II, III and IV, see Conway-Sloane \cite{CoSl}), and Duursma proved the following: 
\begin{thm}[Duursma \cite{Du4}]\label{thm:type4}
All extremal Type IV codes satisfy the Riemann hypothesis. 
\end{thm}

This paper attempts to extend Duursma's theory to other classes of homogeneous polynomials than the weight enumerators of existing codes. Studying carefully the proof of existence of $P(T)$ (see Appendix), we notice that $W_C(x,y)$ in Definition \ref{dfn:zeta} need not be a weight enumerator of an existing code: more essential point is that $W_C(x,y)$ is a homogeneous polynomial of the form 
\begin{equation}\label{eq:homogen}
x^n+\sum_{i=d}^n A_i x^{n-i}y^i \quad (A_i\in {\bf C},\ A_d\ne0).
\end{equation}
In fact, for any polynomial $W(x,y)$ of the form (\ref{eq:homogen}), we can similarly verify existence and uniqueness of $P(T)\in{\bf C}[T]$ such that 
$$\frac{P(T)}{(1-T)(1-qT)}(y(1-T)+xT)^n=\cdots +\frac{W(x,y)-x^n}{q-1}T^{n-d}+ \cdots$$
(the number $q$ should be determined suitably according to what meaning $W(x,y)$ has). This fact is already used in Duursma's papers when he considers the weight enumerators of the MDS codes, but we would like to go even further. 

As a class of homogeneous polynomials, we consider so-called ``formal weight enumerators''. The notion of formal weight enumerators was first introduced by Ozeki \cite{Oz}, in which he deduced a remarkable result in the theory of modular forms, the construction of the Eisenstein series $E_6(z)$ using an example of the formal weight enumerators and the Brou\'e-Enguehard map. The formal weight enumerator $W(x,y)$ resembles the weight enumerators of Type II codes, but is distinguished from them by the formula
$$ W\left(\frac{x+y}{\sqrt 2}, \frac{x-y}{\sqrt 2}\right)=-W(x,y)$$
(see Definition \ref{dfn:fwe}. When $W(x,y)$ is a weight enumerator of a Type II code, we have $W((x+y)/\sqrt 2, (x-y)/\sqrt 2)=W(x,y)$). We consider $W(x,y)$ a weight enumerator of a virtual binary self dual-code, so we set $q=2$. Then we can determine its zeta polynomial $P(T)$. It turns out from the above formula that the functional equation of $P(T)$ becomes
\begin{equation}\label{eq:feq_fwe}
P(T)=-P\Bigl(\frac{1}{2T}\Bigr)2^g T^{2g}
\end{equation}
($g=\frac{n}{2}+1-d$, see Theorem \ref{thm:feq_fwe}). The set of all weight enumerators of Type II codes and all formal weight enumerators forms the invariant polynomial ring ${\bf C}[x,y]^{G_8}$ (see (\ref{eq:group8})), so it follows that we found two different functional equations for the members of ${\bf C}[x,y]^{G_8}$:
\begin{eqnarray*}
P(T)&=&P\Bigl(\frac{1}{2T}\Bigr)2^g T^{2g}\quad (\mbox{Type II codes}),\\
P(T)&=&-P\Bigl(\frac{1}{2T}\Bigr)2^g T^{2g}\quad (\mbox{formal weight enumerators}). 
\end{eqnarray*}
It seems that such a pair of functional equations has never been encountered in the context of number theory. 

What is interesting is that we can formulate an analogue of the Riemann hypothesis for $P(T)$ which satisfy (\ref{eq:feq_fwe}) in a similar way to  the case of the original Duursma theory (see Proposition \ref{prop:roots_of_p(t)} and Definition \ref{dfn:rh_for_fwe}). Moreover there seems to be a structure of the zeta functions for formal weight enumerators quite similar to that of the zeta functions for Type II codes: we can define the extremal property for formal weight enumerators and may guess some relation  between extremal formal weight enumerators and the Riemann hypothesis (see Section 4). 

These results suggest that we should not restrict ourselves to existing linear codes when we consider ``zeta functions for {\it linear codes}'', and that we should take into consideration various other classes of invariant polynomials such as formal weight enumerators. 

For closer observations of formal weight enumerators, we proceed to quantitative investigation for the extremal property of them. For Type II codes, the best possible bound for the minimum distance $d$ is known (see also MacWilliams-Sloane \cite[pp.624-628]{MaSl}): 
\begin{thm}[Mallows-Sloane \cite{MalSl}]\label{thm:mallows-sloane}
For any Type II code of length $n$ and minimum distance $d$, 
$$d\leq 4\left[\frac{n}{24}\right]+4.$$
\end{thm}
The equality holds in the above theorem for an extremal Type II code. Writing a formal weight enumerator $W(x,y)$ in the form (\ref{eq:homogen}), we call $W(x,y)$ extremal if $d$ is the largest possible for given $n$ (see Definition \ref{dfn:ext_fwe}). We deduce the best possible bound for formal weight enumerators analogous to Theorem \ref{thm:mallows-sloane}, which is one of the main results of this paper: 
\begin{thm}\label{thm:bound_fwe}
For any formal weight enumerator $W(x,y)$ of the form (\ref{eq:homogen}), we have
$$d\leq 4\left[\frac{n-12}{24}\right]+4.$$
The equality holds when $W(x,y)$ is extremal. 
\end{thm}
We can see from the discussion of Duursma \cite{Du4} that Theorem \ref{thm:mallows-sloane} is also valid for formal weight enumerators, but it is {\it not the best possible}. So Theorem \ref{thm:bound_fwe} improves Theorem \ref{thm:mallows-sloane} for the case of formal weight enumerators. 

In Section 2, we give a precise definition of formal weight enumerators and give some basic properties of them. In Section 3, we deduce the functional equation of the zeta functions of formal weight enumerators and formulate an analogue of the Riemann hypothesis. Section 4 is devoted to the extremal property of the formal weight enumerators. We observe the relation between the Riemann hypothesis and the extremal property. A proof of Theorem \ref{thm:bound_fwe} is given in Section 5. In Section 6, we provide some numerical data. All the extremal formal weight enumerators of degrees up to 200 are determined, and we mention the Riemann hypothesis again. In Appendix, we give an elementary proof of existence and uniqueness of the zeta polynomial $P(T)$. This fact was first established in Duursma \cite{Du1}, but a detailed proof was omitted. We give an alternative proof, which is also valid for $W(x,y)$ with complex coefficients. 

For part of the results in this paper, see also \cite{Ch1}. The ternary case is dealt with in \cite{Ch2}. 
\section{Formal weight enumerators}
First we introduce the notion of formal weight enumerators: 
\begin{dfn}\label{dfn:fwe}
We call $W(x,y)=\sum_{i=0}^n A_i x^{n-i}y^i \in {\bf C}[x,y]$ a formal weight enumerator if the following conditions are satisfied: 

\smallskip
\noindent(i) $A_i\ne 0$ $\Rightarrow$ $4|i$. \\
(ii) $\displaystyle W\left(\frac{x+y}{\sqrt 2}, \frac{x-y}{\sqrt 2}\right)=-W(x,y)$. 
\end{dfn}
The formal weight enumerators belong to the invariant polynomial ring ${\bf C}[x,y]^{G_8}$ where $G_8$ is defined by
\begin{equation}\label{eq:group8}
G_8:=\left\langle 
\sigma_1=\frac{1-i}{2}
\left(\begin{array}{rr}1 & -1 \\ 1 & 1 \end{array}\right),\ 
\sigma_2=\left(\begin{array}{rr}-i & 0 \\ 0 & 1 \end{array}\right) \right\rangle
\end{equation}
(see Shephard-Todd \cite{ShTo}). The action of a linear transformation 
$\sigma=\left(\begin{array}{cc}a & b \\ c & d \end{array}\right)$ 
to elements $A(x,y)\in{\bf C}[x,y]$ is defined by
$$A(x,y)^\sigma=A(ax+by,cx+dy).$$
It is known that ${\bf C}[x,y]^{G_8}$ is generated by the following two polynomials: 
\begin{eqnarray}
W_8(x,y)&=&x^{8}+14x^4y^4+y^8,\label{eq:ext_hamming}\\
W_{12}(x,y)&=&x^{12}-33x^8y^4-33x^4y^8+y^{12}\label{eq:gen_12}.
\end{eqnarray}
The polynomial $W_8(x,y)$ is the weight enumerator of the extended Hamming code, a well-known example of Type II codes, and $W_{12}(x,y)$ satisfies the condition (ii) of Definition \ref{dfn:fwe}. The condition (i) can be understood by invariance under the action of

$$\sigma_2(\sigma_1^2 \sigma_2^3)^2
=\left(\begin{array}{rr}1 & 0 \\ 0 & i \end{array}\right).$$
Ozeki \cite{Oz} uses $W_{12}(x,y)$ and the Brou\'e-Enguehard map to construct the Eisenstein series $E_6(z)$. 

We can see that general forms of formal weight enumerators are
\begin{equation}\label{eq:fwe_basic}
W_8(x,y)^s W_{12}(x,y)^{2t+1}\quad(s,t\geq 0)
\end{equation}
and their suitable linear combinations. 

The next lemma plays a crucial role in Section 5: 
\begin{lem}\label{lem:property_fwe}
Let $W(x,y)$ be a formal weight enumerator. Then we have the following: \\
(i) $\deg W\equiv 4$ {\rm (mod 8)} and $W(x,y)$ consists of even number of terms. \\
(ii) $W(y,x)=W(x,y)$. 
\end{lem}
\prf Both assertions are straightforward if we assume (\ref{eq:fwe_basic}), but here we show them directly using the action of $\sigma_1$ and $\sigma_2$. \\
(i) Let $\deg W=n$. Then from Definition \ref{dfn:fwe} (ii), 
\begin{eqnarray*}
W(x,y)^{\sigma_1}&=&\left(\frac{1-i}{2}\right)^nW(x-y,x+y)\\
     &=&-\left(\frac{1-i}{2}\right)^nW\left(\frac{(x-y)+(x+y)}{\sqrt 2},\frac{(x-y)-(x+y)}{\sqrt 2}\right)\\
     &=&-\left(\frac{1-i}{\sqrt 2}\right)^n W(x,-y)\\
     &=&(-1)^{n/4+1}W(x,y).
\end{eqnarray*}
In the last equality, we use $W(x,-y)=W(x,y)$, invariance under $\{\sigma_2(\sigma_1^2 \sigma_2^3)^2\}^2$. So if $W(x,y)$ is invariant under $\sigma_1$, we have $2|(\frac{n}{4}+1)$, which leads the congruence. Moreover, the variable $y$ always appears in the form $y^4$ (Definition \ref{dfn:fwe} (i)), so the number of terms of $W(x,y)$ is even. 

\smallskip
\noindent (ii) We can see the formula by the fact that $W(x,y)$ is invariant under 
$$\sigma_2\sigma_1^2 \sigma_2^3=\left(\begin{array}{rr}0 & 1 \\ 1 & 0 \end{array}\right).$$
\hfill\qed
\section{Zeta functions and an analogue of the Riemann hypothesis}
Let $W(x,y)$ be a formal weight enumerator and we write it in the form (\ref{eq:homogen}). We consider $W(x,y)$ a weight enumerator of a virtual binary self dual-code, so we set $q=2$. Then we can determine its zeta polynomial $P(T)$. The first result is the functional equation of $P(T)$: 
\begin{thm}\label{thm:feq_fwe}
The zeta polynomial $P(T)$ of $W(x,y)$ is of degree $2g$ ($g=\frac{n}{2}+1-d$) and satisfies 
\begin{equation}\label{eq:feq_fwe_in_thm}
P(T)=-P\Bigl(\frac{1}{2T}\Bigr)2^g T^{2g}.
\end{equation}
\end{thm}
\prf The fact $\deg P=2g$ can be shown similarly to \cite[p.59]{Du2}. We put 
$$W^\perp(x,y)= W\left(\frac{x+y}{\sqrt 2}, \frac{x-y}{\sqrt 2}\right)$$
(it is sometimes called ``the MacWilliams transform'' of $W(x,y)$). Let $P^\perp(T)$ be the zeta polynomial of $W^\perp(x,y)$. Then we can show 
$$P^\perp(T)=P\left(\frac{1}{2T}\right)2^gT^{2g} \quad(g=\textstyle\frac{n}{2}+1-d).$$
The proof is similar to that of \cite[p.59]{Du2}. The polynomial $P^\perp(T)$ must coincides with the zeta polynomial of $-W(x,y)$ by Definition \ref{dfn:fwe} (ii), so we get the desired formula. \qed

\medskip
Next we examine the distribution of the roots of $P(T)$: 
\begin{prop}\label{prop:roots_of_p(t)}
Let $P(T)$ be the zeta polynomial of a formal weight enumerator. Then we can arrange the roots of $P(T)$ as follows: 
$$\alpha_1, \frac{1}{2\alpha_1}, \cdots, \alpha_s, \frac{1}{2\alpha_s}, 
\frac{1}{\sqrt 2}, \cdots, \frac{1}{\sqrt 2}, -\frac{1}{\sqrt 2}, \cdots, -\frac{1}{\sqrt 2}$$
for some $s\in{\bf N}$, $\alpha_j\ne \pm 1/\sqrt 2$ ($1\leq j\leq s$), both $1/\sqrt 2$ and $-1/\sqrt 2$ occur in odd multiplicities. 
\end{prop}
\prf From Theorem \ref{thm:feq_fwe}, $P(\alpha)=0$ is equivalent to $P(1/2\alpha)=0$, so the roots $\alpha$ and $1/2\alpha$ appear in pairs. If $\alpha=
1/2\alpha$, then $\alpha=\pm1/\sqrt 2$, so we can arrange the roots of $P(T)$ as follows: 
$$\alpha_1, \frac{1}{2\alpha_1}, \cdots, \alpha_s, \frac{1}{2\alpha_s}, 
\frac{1}{\sqrt 2}, \cdots, \frac{1}{\sqrt 2}, -\frac{1}{\sqrt 2}, \cdots, -\frac{1}{\sqrt 2}$$
for some $s\in{\bf N}$, $\alpha_j\ne \pm 1/\sqrt 2$ ($1\leq j\leq s$). Suppose the multiplicity of $1/\sqrt 2$ equals $m$ and that of $-1/\sqrt 2$ equals $n$. Let
$$P(T)=a_0+a_1T+\cdots +a_{2g}T^{2g}.$$
Then we have
$$P\left(\frac{1}{2T}\right)2^gT^{2g}=2^{-g}a_{2g}+2^{-g+1}a_{2g-1}T+\cdots +2^g a_0 T^{2g}.$$
Comparing the coefficients of both sides of (\ref{eq:feq_fwe_in_thm}), we get $2^g a_0=-a_{2g}$ and 
$$P(T)=a_{2g}\left(T^{2g}+\cdots +\frac{a_1}{a_{2g}}T-\frac{1}{2^g}\right).$$
Therefore the product of all roots of $P(T)$ equals $(-1)^{2g}(-1/2^g)=-1/2^g$, and we have
$$\frac{(-1)^n}{2^{s+m/2+n/2}}=\frac{-1}{2^g}.$$
We can see from this formula that both $m$ and $n$ are odd. \qed

\medskip
\noindent\rem In the cases of the zeta polynomials for algebraic curves or existing self-dual codes (over $\ff q$), the multiplicities of $\pm 1/\sqrt  q$ are even. It is one of the different points of them from the formal weight enumerators. 

\medskip
\noindent Taking Proposition \ref{prop:roots_of_p(t)} into account, we see that it is appropriate to formulate the Riemann hypothesis for formal weight enumerators as follows: 
\begin{dfn}\label{dfn:rh_for_fwe}
A formal weight enumerator $W(x,y)$ satisfies the Riemann hypothesis if all the zeros of its zeta polynomial $P(T)$ have the same absolute value $1/\sqrt 2$. 
\end{dfn}

The Riemann hypothesis for existing self-dual codes is introduced in Definition \ref{dfn:RH} (see also \cite[Definition 4.1]{Du3}). Duursma observes that good codes tend to satisfy the Riemann hypothesis (see \cite[Section 1]{Du2} and \cite[Abstract]{Du3}). Definition \ref{dfn:rh_for_fwe} is just a formal analogy to Definition \ref{dfn:RH}, but we nevertheless find the evidence that the Riemann hypothesis for formal weight enumerators also implies some ``good'' property of them --- the extremal property. 
\section{Extremal formal weight enumerators}
In this section, we introduce the notion of extremal formal weight enumerators and observe the relation between them and the Riemann hypothesis. We write formal weight enumerators in the form (\ref{eq:homogen}). 
\begin{dfn}\label{dfn:ext_fwe}
We call a formal weight enumerator $W(x,y)$ ($\deg W=n$) ``extremal'' if $d$ is the largest among all formal weight enumerators of degree $n$. 
\end{dfn}
From the discussion of the ring ${\bf C}[x,y]^{G_8}$ (see (\ref{eq:fwe_basic}) in particular), when $\deg W\leq 28$, $W_{12}(x,y)$, $W_8(x,y) W_{12}(x,y)$ and $W_8(x,y)^2 W_{12}(x,y)$ are themselves extremal, but when $\deg W\geq 36$, there always exist at least two different formal weight enumerators, so we can eliminate the terms with small powers in $y$ and construct the extremal formal weight enumerator. The process is illustrated as follows: 
\begin{exam}\label{exam:deg36}\rm 
$\deg W=36$. There are two formal weight enumerators of the form (\ref{eq:fwe_basic}): 
\begin{eqnarray*}
W_8(x,y)^3 W_{12}(x,y)&=&x^{36}+9x^{32}y^4-828x^{28}y^8-\cdots\\
W_{12}(x,y)^3&=&x^{36}-99x^{32}y^4+3168x^{28}y^8-\cdots.
\end{eqnarray*}
In this case, 
$$\frac{11}{12}W_8(x,y)^3 W_{12}(x,y)+\frac{1}{12}W_{12}(x,y)^3=x^{36}-495x^{28}y^8-19005x^{24}y^{12}-\cdots$$
is extremal. 
\end{exam}
There seems to be a strong resemblance between the structure of the zeta functions for formal weight enumerators and that of the zeta functions for Type II codes. In fact, the following examples suggest that the relation between the extremal property and the Riemann hypothesis is the same as the case of  Type II codes: 
\begin{exam}\label{exam:rh-1}\rm 
Numerical examples of the Riemann hypothesis for formal weight enumerators of small degrees. In the column ``RH'' ($=$ the Riemann hypothesis), ``T'' ($=$ true) means that the Riemann hypothesis holds for $W(x,y)$, and ``F'' ($=$ false) means that it does not hold. In the examples below, $W(x,y)$ with $\deg W=12, 20, 28$ are extremal, but others are not (for all $W(x,y)$, the number $d$ in (\ref{eq:homogen}) is equal to 4). 

\medskip
\begin{center}
\begin{tabular}{|l|c|c|}\hline
$W(x,y)$ & $\deg W$ & RH\\\hline
$W_{12}(x,y)$ & 12 & T\\
$W_8(x,y) W_{12}(x,y)$ & 20 & T\\
$W_8(x,y)^2 W_{12}(x,y)$ & 28 & T\\
$W_8(x,y)^3 W_{12}(x,y)$ & 36 & F\\
$W_{12}(x,y)^3$ & 36  & F\\
$W_8(x,y)^4 W_{12}(x,y)$ & 44 & F\\
$W_8(x,y) W_{12}(x,y)^3$ & 44 & F\\
$W_8(x,y)^5 W_{12}(x,y)$ & 52 & F\\
$W_8(x,y)^2 W_{12}(x,y)^3$ & 52 & F\\
$W_8(x,y)^6 W_{12}(x,y)$ & 60 & F\\
$W_8(x,y)^3 W_{12}(x,y)^3$ & 60 & F\\
$W_{12}(x,y)^5$ & 60 & F\\\hline
\end{tabular}
\end{center}
We give the zeta polynomials $P_{12}(T)$, $P_{20}(T)$ and $P_{28}(T)$ for $W_{12}(x,y)$, $W_8(x,y) W_{12}(x,y)$ and $W_8(x,y)^2 W_{12}(x,y)$, respectively: 
\begin{eqnarray*}
P_{12}(T)&=&\frac{1}{15}(2T^2-1)(2T^2+1)(2T^2+2T+1),\\
P_{20}(T)&=&\frac{1}{255}(2T^2-1)(2T^2+2T+1)(2T^2+1)(16T^8+1),\\
P_{28}(T)&=&\frac{1}{4095}(2T^2-1)(2T^2+2T+1)(2T^2+1)(4T^4-2T^2+1)(4T^4+2T^2+1)\\
 & &\cdot(4T^4+4T^3+2T^2+2T+1)(4T^4-4T^3+2T^2-2T+1).
\end{eqnarray*}
\end{exam}
When $\deg W\geq 36$, we can construct the extremal formal weight enumerators in the same manner as in Example \ref{exam:deg36}: 
\begin{exam}\label{exam:rh-2}\rm 
Numerical examples of the Riemann hypothesis for some extremal formal weight enumerators. 

\begin{center}
\begin{tabular}{|c|c|c|c|}\hline
$W(x,y)$ & $\deg W$ & $d$ & RH\\\hline
$\frac{11}{12} W_8(x,y)^3 W_{12}(x,y) + \frac{1}{12} W_{12}(x,y)^3$ & 36 & 8 & T\\
$\frac{85}{108} W_8(x,y)^4 W_{12}(x,y) + \frac{23}{108} W_8(x,y) W_{12}(x,y)^3$ 
& 44 & 8 & T\\
$\frac{71}{108} W_8(x,y)^5 W_{12}(x,y) + \frac{37}{108} W_8(x,y)^2 W_{12}(x,y)^3$ 
& 52 & 8 & T\\
$\frac{1045}{1944} W_8(x,y)^6 W_{12}(x,y) 
+ \frac{880}{1944} W_8(x,y)^3 W_{12}(x,y)^3
+ \frac{19}{1944} W_{12}(x,y)^5$ 
& 60 & 12 & T\\\hline
\end{tabular}
\end{center}
The zeta polynomial $P_{36}(T)$ for $\frac{11}{12} W_8(x,y)^3 W_{12}(x,y) + \frac{1}{12} W_{12}(x,y)^3$ is
\begin{eqnarray*}
P_{36}(T)&=&{\frac {1}{11920740}}\,\left (2\,{T}^{2}-1\right ) (199680\,{T}^{
20}+599040\,{T}^{19}+1098240\,{T}^{18}+1497600\,{T}^{17}\\
& &+1683904\,{T}^{16}+1630400\,{T}^{15}+1410176\,{T}^{14}+1116384\,{T}^{13}\\
& &+832384\,{T}^{12}+598544\,{T}^{11}+424720\,{T}^{10}+299272\,{T}^{9}+208096\,{T}^{8}+139548\,{T}^{7}\\
& &+88136\,{T}^{6}+50950\,{T}^{5}+26311\,{T}^{4}+11700\,{T}^{3}+4290\,{T}^{2}+1170\,T+195 )
\end{eqnarray*}
(we omit the zeta polynomials for others). 
\end{exam}
Seeing these examples, we may ask a question of whether extremal formal weight enumerators satisfy the Riemann hypothesis. 

For more numerical data, the reader is referred to Section 6. 
\section{Proof of Theorem \ref{thm:bound_fwe}}
Our proof of Theorem \ref{thm:bound_fwe} is an application of the technique developed in \cite[Section 2]{Du4}. First we introduce some notations and a lemma from \cite{Du4}. 

For a linear transformation $\sigma=
\left(\begin{array}{cc} a & b \\ c & d \end{array}\right)$ we define the relation between two pairs of variables $(x,y)$ and $(u,v)$ as follows in accordance with \cite[Section 2]{Du4}: 
\begin{equation}\label{eq:transform_uv_xy}
(u,v)=(x,y)\sigma=(ax+cy, bx+dy)
\end{equation}
(note that the action of $\sigma$ is different from that of Section 2). We introduce the matching transformation of differential operators
$$\left(\frac{\partial}{\partial x}, \frac{\partial}{\partial y}\right)
=\left(\frac{\partial}{\partial u}, \frac{\partial}{\partial v}\right)
\sigma^{{\rm T}}
=\left(a\frac{\partial}{\partial u}+b\frac{\partial}{\partial v},
 c\frac{\partial}{\partial u}+d\frac{\partial}{\partial v}\right),$$
where $\sigma^{{\rm T}}$ means the transposed matrix of $\sigma$. 

Let $a(x,y)$, $p(x,y)$ and $A(x,y)$ be arbitrary homogeneous polynomials over ${\bf C}$. We denote the differential operator $p(\partial/\partial x, \partial/\partial y)$ by $p(x,y)(D)$. Then we have
\begin{lem}\label{lem:duursma-1}
$$p((u,v)\sigma^{{\rm T}})(D)A(u,v)=p(x,y)(D)A((x,y)\sigma).$$
\end{lem}
\prf This is Lemma 1 of Duursma \cite{Du4}. \qed

\medskip
The basic idea is that we would like to find a relation of the form
\begin{equation}\label{eq:basic_idea}
a(x,y)|p(x,y)(D)A(x,y)
\end{equation}
between a (formal) weight enumerator $A(x,y)$ and some polynomials $a(x,y)$, $p(x,y)$. Then we can say $\deg a$ is less than or equal to the degree of the right hand side. The degrees of the terms in (\ref{eq:basic_idea}) contain parameters such as the code length $n$ and the minimum distance $d$. If we can find good $a(x,y)$ and $p(x,y)$, then the inequality of the degrees becomes straightforwardly a bound of $d$ in terms of $n$. By this method, Duursma obtains an alternative proof of the Mallows-Sloane bounds for Types I through IV (see Theorem 3 and Section 1.1 of \cite{Du4}). We apply this method to formal weight enumerators. 

Let $W(x,y)$ be a formal weight enumerator of degree $n$. The number $d$ is the same as in (\ref{eq:homogen}). Then we have the following two propositions: 
\begin{prop}\label{prop:basic-1}If $d\geq 8$, 
$$(xy)^{d-5}(x^4-y^4)^{d-5}|xy(x^4-y^4)(D)W(x,y).$$
\end{prop}
\prf It can be shown similarly to \cite[Lemma 2]{Du4} (recall Definition \ref{dfn:fwe} (i) and Lemma \ref{lem:property_fwe} (i)). \qed
\begin{prop}\label{prop:basic-2}If $d\geq 8$, 
$$(x^4+y^4)(x^4+6x^2y^2+y^4)|xy(x^4-y^4)(D)W(x,y).$$
\end{prop}
\prf By Definition \ref{dfn:fwe} (i) and Lemma \ref{lem:property_fwe}, $W(x,y)$ can be written in the form
\begin{equation}\label{eq:fwe_symm}
W(x,y)=x^n+y^n+\sum_{j=d/4}^{\frac{n-4}{8}} A_{4j}(x^{n-4j}y^{4j}+x^{4j}y^{n-4j}).
\end{equation}
Then we can easily verify 
\begin{eqnarray*}
xy(x^4-y^4)(D)W(x,y)&=&
\sum_{j=d/4}^{\frac{n-4}{8}}A_{4j}
\{ (n-4j)_4 \cdot 4j(n-4j-4)(x^{n-4j-5}y^{4j-1}-x^{4j-1}y^{n-4j-5})\\
& &+(4j)_4 (n-4j)\cdot 4(j-1)(x^{4j-5}y^{n-4j-1}-x^{n-4j-1}y^{4j-5})\},
\end{eqnarray*}
where $(a)_n=a(a-1)\cdots (a-n+1)$. We can put $n=8s+4$ ($s\in{\bf N}$) by Lemma \ref{lem:property_fwe} (i) and have
$$x^{n-4j-5}y^{4j-1}-x^{4j-1}y^{n-4j-5}=(xy)^{4j-1}(x^{8(s-j)}-y^{8(s-j)})$$
which is divisible by $x^4+y^4$. Similarly, 
$$(x^4+y^4)|(x^{4j-5}y^{n-4j-1}-x^{n-4j-1}y^{4j-5}).$$
Next we apply Lemma \ref{lem:duursma-1} with $A(x,y)=W(x,y)$, 
$\sigma=\left(\begin{array}{rr}1 & 1 \\ 1 & -1 \end{array}\right)$ and $p(x,y)=xy(x^4-y^4)$. By Definition \ref{dfn:fwe} (ii), we have
$$p(x,y)(D)W((x,y)\sigma)=p(x,y)(D)W(x+y, x-y)=
-(\sqrt 2)^n p(x,y)(D)W(x,y).$$
Since $(x^4+y^4)|p(x,y)(D)W(x,y)$, $p((u,v)\sigma^{{\rm T}})(D)W(u,v)$ is divisible by the polynomial 
$$\frac{1}{8}(u^4+6u^2v^2+v^4),$$
the image of $x^4+y^4$ by the transformation $\sigma$. On the other hand, we have 
$$p((u,v)\sigma^{{\rm T}})=p(u+v, u-v)=8uv(u^4-v^4)=8p(u,v).$$
Therefore
$$(x^4+6x^2y^2+y^4)|xy(x^4-y^4)(D)W(x,y).$$
Two polynomials $x^4+y^4$ and $x^4+6x^2y^2+y^4$ are coprime, so we get the desired result. \qed

\medskip
\noindent Propositions \ref{prop:basic-1} and \ref{prop:basic-2} bring the following: 
\begin{thm}\label{thm:divisibility}
For any formal weight enumerator with $d\geq 8$, we have
\begin{equation}\label{eq:divisibility}
(xy)^{d-5}(x^4-y^4)^{d-5}(x^4+y^4)(x^4+6x^2y^2+y^4)|xy(x^4-y^4)(D)W(x,y).
\end{equation}
\end{thm}
Theorem \ref{thm:bound_fwe} is deduced from the above theorem. Comparing the degrees in both sides of (\ref{eq:divisibility}), we have
$$6(d-5)+8\leq n-6, \quad 6d\leq n-8+24.$$
Since $4|d$, we put $d=4d'$. Then, $d'\leq (n-8)/24+1$. Since $d'\in{\bf Z}$, this is equivalent to
$$d'\leq \left[\frac{n-8}{24}\right]+1, \quad d\leq 4\left[\frac{n-8}{24}\right]+4.$$
For $n\equiv 4$ (mod 8), we can easily verify $[(n-8)/24]=[(n-12)/24]$. Hence we get the inequality in Theorem \ref{thm:bound_fwe}. 

Finally, we consider extremal formal weight enumerators. Using two generators (\ref{eq:ext_hamming}) and (\ref{eq:gen_12}), we can write a formal weight enumerator in the form $\sum a_s W_8(x,y)^s W_{12}(x,y)^{2t+1}$ ($8s+12(2t+1)=n$), but noting that
$$\frac{1}{108}\{W_8(x,y)^3-W_{12}(x,y)^2\}=x^4y^4(x^4-y^4)^4,$$
we use $W'_{24}(x,y):=x^4y^4(x^4-y^4)^4$ instead of $W_8(x,y)^3$. We classify formal weight enumerators into the following three sequences:

\medskip
(I) $\displaystyle W(x,y)=
    \sum_{r=0}^m a_r W'_{24}(x,y)^r W_{12}(x,y)^{2m-2r+1}$ 
    ($m\geq 0$, $\deg W=24m+12$)

\medskip
(II) $\displaystyle W(x,y)=
    \sum_{r=0}^m a_r W'_{24}(x,y)^r \cdot W_8(x,y)W_{12}(x,y)^{2m-2r+1}$ 
    ($m\geq 0$, $\deg W=24m+20$)

\medskip
(III) $\displaystyle W(x,y)=
    \sum_{r=0}^m a_r W'_{24}(x,y)^r \cdot W_8(x,y)^2W_{12}(x,y)^{2m-2r+1}$ 
    ($m\geq 0$, $\deg W=24m+28$)

\medskip
\noindent In each sequence, the number of summands is $m+1$, so we can eliminate the terms with $y^4$, $y^8$, $\cdots$, $y^{4m}$, choosing suitable $a_r$'s. Thus an extremal formal weight enumerator $W(x,y)$ must satisfy $d\geq 4m+4$. On the other hand, for $n=24m+12$, $24m+20$ and $24m+28$, the formula $4[(n-12)/24]+4$ gives the same value $4m+4$. Therefore the equality holds in Theorem \ref{thm:bound_fwe} if and only if $W(x,y)$ is extremal. \qed

\medskip
\noindent\rem The inequality of Theorem \ref{thm:bound_fwe} can also be proved using the theory of modular forms (Siegel's theorem, see Corollary 1 on p.134 of \cite{Hi}). The author expresses his gratitude to Professor Eiichi Bannai for pointing out this viewpoint. Here we showed it directly, using less prerequisites. 
\section{Numerical examples}
In this section, we write formal weight enumerators in the form
$$W(x,y)=x^n+y^n+\sum_{j=d/4}^{\frac{n-4}{8}} A_{4j}(x^{n-4j}y^{4j}+x^{4j}y^{n-4j})$$
(see (\ref{eq:fwe_symm})). We list all extremal formal weight enumerators with degrees up to 200. One can verify numerically that all these $W(x,y)$ satisfy the Riemann hypothesis.

\medskip
\noindent\begin{tabular}{r|r|l}
$n$ & $d$ & $A_j$\\\hline
12 & 4 & $A_4=-33$\\

20 & 4 & $A_4=-19$, $A_8=-494$\\

28 & 4 & $A_4=-5$, $A_8=-759$, $A_{12}=-7429$\\

36 & 8 & $A_8=-495$, $A_{12}=-19005$, $A_{16}=-111573$\\

44 & 8 & $A_8=-172$, $A_{12}=-20167$, $A_{16}=-397234$, $A_{20}=-1679580$\\

52 & 8 & $A_8=-45$, $A_{12}=-12313$, $A_{16}=-617766$, $A_{20}=-7509579$, 
         $A_{24}=-25414730$\\

60 & 12 & $A_{12}=-5605$, $A_{16}=-554895$, $A_{20}=-15622728$, 
          $A_{24}=-15622728$, \\
   &    & $A_{28}=-386391885$ \\

68 & 12 & $A_{12}=-1742$, $A_{16}=-342705$, $A_{20}=-19114095$, 
          $A_{24}=-350115870$,\\
   &    & $A_{28}=-2321175604$, $A_{32}=-5899184577$\\

76 & 12 & $A_{12}=-455$, $A_{16}=-157788$, $A_{20}=-15862224$, 
          $A_{24}=-548257680$,\\
   &    & $A_{28}=-7249325900$, $A_{32}=-39225987090$, 
          $A_{36}=-90399362336$\\

84 & 16 & $A_{16}=-62748$, $A_{20}=-9720711$, $A_{24}=-584058384$, 
          $A_{28}=-13904530512$, \\
   &    & $A_{32}=-142025799762$, $A_{36}=-652738809996$, 
          $A_{40}=-1389760273440$\\

92 & 16 & $A_{16}=-18564$, $A_{20}=-4750200$, $A_{24}=-459415320$,
          $A_{28}=-18276716575$,\\
   &    & $A_{32}=-322921346202$, $A_{36}=-2673263269352$,
          $A_{40}=-10743644373000$, \\
   &    & $A_{44}=-21425802199620$\\

100 & 16 & \small{$A_{16}=-4845$, $A_{20}=-1906569$, $A_{24}=-283417680$, 
           $A_{28}=-17759231600$,}\\
    &    & \small{$A_{32}=-508083348150$, $A_{36}=-7024438840315$, 
           $A_{40}=-48836439768144$,}\\
    &    & \small{$A_{44}=-175426727301360$, $A_{48}=-331136219602650$}\\

108 & 20 & \small{$A_{20}=-707805$, $A_{24}=-142345845$, 
           $A_{28}=-13483875924$, $A_{32}=-592827012162$,}\\
    &    & \small{$A_{36}=-12912473819940$, $A_{40}=-145356811396020$, 
           $A_{44}=-872012120596695$,} \\
    &    & \small{$A_{48}=-2847442918238100$, $A_{52}=-5128868476748502$} \\

116 & 20 & \small{$A_{20}=-203665$, $A_{24}=-61035220$, $A_{28}=-8325627860$,           $A_{32}=-539970141755$,}\\
    &    & \small{$A_{36}=-17686996366552$, $A_{40}=-306089064693104$, 
           $A_{44}=-2892828876388720$,}\\
    &    & \small{$A_{48}=-15295065179581320$, 
           $A_{52}=-46010051597767125$,}\\
    &    & \small{$A_{56}=-79592918004050552$}\\

124 & 20 & \small{$A_{20}=-53130$, $A_{24}=-22587526$, 
           $A_{28}=-4320543557$, $A_{32}=-399142595247$,}\\
    &    & \small{$A_{36}=-18933376280702$, $A_{40}=-483099940106322$, 
           $A_{44}=-6866557381558152$,}\\
    &    & \small{$A_{48}=-55812566735037112$, 
           $A_{52}=-264497659587312580$,}\\
    &    & \small{$A_{56}=-740865653481227100$, 
           $A_{60}=-1237298135226392525$}\\

132 & 24 & \small{$A_{24}=-8055190$, $A_{28}=-1920233370$, 
           $A_{32}=-246940710159$,}\\
    &    & \small{$A_{36}=-16432390247845$, $A_{40}=-597323869611810$, 
           $A_{44}=-12294641192051790$,}\\
    &    & \small{$A_{48}=-147341787206768440$, 
           $A_{52}=-1050128724293493768$,} \\
    &    & \small{$A_{56}=-4521938050758467196$, 
           $A_{60}=-11897286760705846500$,}\\
    &    & \small{$A_{64}=-19263884178133617165$}\\

140 & 24 & \small{$A_{24}=-2276820$, $A_{28}=-757465210$, 
           $A_{32}=-130844977725$,}\\
    &    & \small{$A_{36}=-11905432407540$, $A_{40}=-598991736062244$, 
           $A_{44}=-17302162719073995$,}\\
    &    & \small{$A_{48}=-295390461241721370$, 
           $A_{52}=-3048727593848694660$,}\\
    &    & \small{$A_{56}=-19356087295370216360$, 
           $A_{60}=-76589231536103045004$,}\\
    &    & \small{$A_{64}=-190646163336241077465$, 
           $A_{68}=-300342296944408633320$}\\
\end{tabular}

\noindent\begin{tabular}{r|r|l}
$n$ & $d$ & $A_j$\\\hline
148 & 24 & \small{$A_{24}=-593775$, $A_{28}=-265687763$, 
           $A_{32}=-60626687121$, $A_{36}=-7369909114334$,}\\
    &    & \small{$A_{40}=-500856000132642$, $A_{44}=-19772475234348690$, 
           $A_{48}=-467244030396146094$,}\\
    &    & \small{$A_{52}=-6767638484266988559$, 
           $A_{56}=-61210009702624439541$,}\\
    &    & \small{$A_{60}=-350734164704203152873$, 
           $A_{64}=-1287227549859336458945$,}\\
    &    & \small{$A_{68}=-3049794439066320485652$, 
           $A_{72}=-4688511639130106191404$}\\

156 & 28 & \footnotesize{$A_{28}=-92385735$, $A_{32}=-24801039900$, 
           $A_{36}=-3974981408320$, $A_{40}=-357147810452160$,}\\
    &    & \footnotesize{$A_{44}=-18838478219043600$, 
           $A_{48}=-601305738440206800$, $A_{52}=-11903058986903663040$,}\\
    &    & \footnotesize{$A_{56}=-149015116832397271872$, 
           $A_{60}=-1198277004601232375820$, }\\
    &    & \footnotesize{$A_{64}=-6264996268248174364050$, 
           $A_{68}=-21495202711204378665600$,}\\
    &    & \footnotesize{$A_{72}=-48721749082300406652800$, 
           $A_{76}=-73273963704290809308576$} \\

164 & 28 & \footnotesize{$A_{28}=-25790512$, $A_{32}=-9233515116$, 
           $A_{36}=-1895408784087$, $A_{40}=-221270131479600$,}\\
    &    & \footnotesize{$A_{44}=-15290974644495840$, 
           $A_{48}=-645538939197188400$, $A_{52}=-17073473360477368176$,}\\
    &    & \footnotesize{$A_{56}=-288709937439922777568$, 
           $A_{60}=-3172914113669269549776$,}\\
    &    & \footnotesize{$A_{64}=-22961929471005530775042$, 
           $A_{68}=-110560597088721242234940$,}\\
    &    & \footnotesize{$A_{72}=-356999855309914228424400$, 
           $A_{76}=-777499897249986858851904$,}\\
    &    & \footnotesize{$A_{80}=-1146350001532072178176992$}\\

172 & 28 & \footnotesize{$A_{28}=-6724520$, $A_{32}=-3117932532$, 
           $A_{36}=-810694410656$, $A_{40}=-120949446496032$,}\\
    &    & \footnotesize{$A_{44}=-10765336143514095$, 
           $A_{48}=-590217737326799400$, $A_{52}=-20452419017792980416$,}\\
    &    & \footnotesize{$A_{56}=-457425536986817800128$, 
           $A_{60}=-6716660868222566243736$,}\\
    &    & \footnotesize{$A_{64}=-65655567144166964191370$, 
           $A_{68}=-432049241644493521436640$,}\\
    &    & \footnotesize{$A_{72}=-1931022707448499094213472$, 
           $A_{76}=-5901887812416818000710932$,}\\
    &    & \footnotesize{$A_{80}=-12396360119130854562120984$, 
           $A_{84}=-17951455639954237535022720$}\\

180 & 32 & \footnotesize{$A_{32}=-1066449780$, $A_{36}=-312308624680$, 
           $A_{40}=-59122804799520$,}\\
    &    & \footnotesize{$A_{44}=-6673483559682720$, 
           $A_{48}=-467612005564042920$, $A_{52}=-20870036989598966895$,}\\
    &    & \footnotesize{$A_{56}=-606173904040335768000$, 
           $A_{60}=-11661733376109411334080$,}\\
    &    & \footnotesize{$A_{64}=-150770765081343953853450$, 
           $A_{68}=-1325608423224821146949400$,}\\
    &    & \footnotesize{$A_{72}=-8002717084203811002888800$, 
           $A_{76}=-33428025572996839645697760$,}\\
    &    & \footnotesize{$A_{80}=-97185266328167090509747608$, 
           $A_{84}=-197504585877104399583402900$,}\\
    &    & \footnotesize{$A_{88}=-281360756340249766957353600$}\\

188 & 32 & \footnotesize{$A_{32}=-294998484$, $A_{36}=-111445052576$, $A_{40}=-26099681914240$,}\\
& & \footnotesize{$A_{44}=-3689805827584200$, $A_{48}=-325758365331044360$, $A_{52}=-18445863562437861440$,}\\
& &  \footnotesize{$A_{56}=-684711621741553803808$, $A_{60}=-16966016533937895991567$,}\\
& & \footnotesize{$A_{64}=-284860403664169699215850$, $A_{68}=-3281429481419712054419680$,}\\
& & \footnotesize{$A_{72}=-26201199634385650653064000$, $A_{76}=-146228109419292198131352376$,}\\
& & \footnotesize{$A_{80}=-574233656555043347026617464$, $A_{84}=-1594905997813459117303489280$,}\\
& & \footnotesize{$A_{88}=-3144907665495143232420795488$, $A_{92}=-4413459725977141147435788980$}\\

196 & 32 & \footnotesize{$A_{32}=-76904685$, $A_{36}=-36580739273$, 
           $A_{40}=-10518107102912$,}\\
    &    & \footnotesize{$A_{44}=-1839425553426240$, 
           $A_{48}=-202055163118166800$, $A_{52}=-14322874390073948080$,}\\
    &    & \footnotesize{$A_{56}=-669909048673946917056$, 
           $A_{60}=-21059701862665139595584$,}\\
    &    & \footnotesize{$A_{64}=-451896822437175266776950$, 
           $A_{68}=-6704657067863573109515835$,}\\
    &    & \footnotesize{$A_{72}=-69523694441231755530686400$, 
           $A_{76}=-508364721051285300494140992$,}\\
    &    & \footnotesize{$A_{80}=-2640396051179528070347697168$, 
           $A_{84}=-9798327455362797587376378160$,}\\
    &    & \footnotesize{$A_{88}=-26097134996226616755233810880$, 
           $A_{92}=-50053424263353462264963577920$,}\\
    &    & \footnotesize{$A_{96}=-69281975548885761832168515738$}\\\hline
\end{tabular}
\section{Appendix --- an elementary proof of existence of $P(T)$}
Existence and uniqueness of the zeta polynomial $P(T)$ for a linear code was first established in \cite[Section 9]{Du1}, but detailed proof is not given. Here we give an alternative, elementary proof, including the case $W(x,y)\in{\bf C}[x,y]$. 

Suppose $W(x,y)\in{\bf C}[x,y]$ is a polynomial of the form (\ref{eq:homogen}). First note that 

\bigskip
$\displaystyle f(T):=\frac{1}{(1-T)(1-qT)}(y(1-T)+xT)^n$
\begin{eqnarray*}
&=&(1+T+T^2+\cdots)(1+qT+q^2T^2+\cdots)((x-y)T+y)^n\\
&=&(1+c_1 T+c_2 T^2+\cdots)\Bigl\{\sum_{j=0}^n \Bigl({n \atop j}\Bigr) (x-y)^jy^{n-j}T^j\Bigr\}
\end{eqnarray*}
for some $c_j\in{\bf N}$. Expanding the last formula, we find for some integers $b_{ij}$, 
$$\begin{array}{lr}
\mbox{the constant term} & =y^n,\\
\mbox{the coefficient of } T& =nxy^{n-1}+(c_1-n)y^n,\\
\cdots\cdots & \cdots\cdots\cdots\\
\mbox{the coefficient of } T^i & =b_{i0}x^iy^{n-i}+b_{i1}x^{i-1}y^{n-i+1}+\cdots +b_{ii}y^n,\\
\cdots\cdots & \cdots\cdots\cdots\\
\mbox{the coefficient of } T^{n-d} & =b_{n-d,0}x^{n-d}y^{d}+b_{n-d,1}x^{n-d-1}y^{d+1}+\cdots +b_{n-d,n-d}y^n.
\end{array}$$
Let $a_0$, $a_1$, $\cdots$, $a_{n-d}\in{\bf C}$ and we form a function $F(T):=(a_0+a_1 T+\cdots +a_{n-d}T^{n-d})f(T)$. Then the coefficient of $T^{n-d}$ of $F(T)$ is 
\begin{eqnarray}
a_{n-d}y^n& &\nonumber\\
+a_{n-d-1}\{nxy^{n-1}+(c_1-n)y^n\}& &\nonumber\\
\cdots\cdots\cdots\qquad\qquad & &\nonumber\\
+a_i\{b_{i0}x^iy^{n-i}+b_{i1}x^{i-1}y^{n-i+1}+\cdots +b_{ii}y^n\}& &\nonumber\\
\cdots\cdots\cdots\qquad\qquad & &\nonumber\\
+a_0\{b_{n-d,0}x^{n-d}y^{d}+b_{n-d,1}x^{n-d-1}y^{d+1}+\cdots +b_{n-d,n-d}y^n\}.& &
\label{eq:coeff_t^n-d}
\end{eqnarray}
On the other hand, since $(W(x,y)-x^n)/(q-1)=(A_dx^{n-d}y^d+\cdots +A_ny^n)/(q-1)$, we can determine $a_0$, $a_1$, $\cdots$, $a_{n-d}$ so that (\ref{eq:coeff_t^n-d}) coincides with $(W(x,y)-x^n)/(q-1)$ (the system of linear equations for determining $a_0$, $a_1$, $\cdots$, $a_{n-d}$ has a regular coefficient matrix). So we can always determine the zeta polynomial $P(T)$ uniquely as $P(T)=a_0+a_1 T+\cdots +a_{n-d}T^{n-d}$. \qed

\end{document}